\documentclass{elsart3-1}

\usepackage{amssymb}
\usepackage[english,francais]{babel}

\newtheorem{theorem}{Theorem}[section]
\newtheorem{lemma}[theorem]{Lemma}
\newtheorem{e-proposition}[theorem]{Proposition}

\newtheorem{e-definition}[theorem]{Definition\rm}
\newtheorem{remark}{\it Remark\/}

\newtheorem{theoreme}{Th\'eor\`eme}[section]

\newtheorem{proposition}[theoreme]{Proposition}

\renewcommand{\theequation}{\arabic{equation}}
\usepackage{amsmath}
\newcommand{\RR}{\mathbb{R}}

\newcommand{\cT}{\mathcal{T}}
\newcommand{\abs}[1]{\left |#1\right |}
\newcommand{\dist}{\operatorname{dist}}
\newcommand{\norm}[1]{\left \|#1\right \|}

\def\og{\leavevmode\raise.3ex\hbox{$\scriptscriptstyle\langle\!\langle$~}}
\def\fg{\leavevmode\raise.3ex\hbox{~$\!\scriptscriptstyle\,\rangle\!\rangle$}}

\begin{document}

\begin{frontmatter}


\selectlanguage{english}
\title{Boundary singularities of positive solutions of some nonlinear
  elliptic equations} 

\vspace{-2.6cm}

\selectlanguage{francais}
\title{Singularit\'es au bord de solutions positives d'\'equations
  elliptiques non-lin\'eaires} 


\selectlanguage{english}
\author{Marie-Fran\c coise Bidaut-V\'eron}, 
\ead{veronmf@univ-tours.fr}
\author{Augusto C. Ponce},
\ead{ponce@lmpt.univ-tours.fr}
\author{Laurent V\'eron}
\ead{veronl@univ-tours.fr}

\address{Laboratoire de Math\'ematiques et Physique Th\'eorique, CNRS
UMR 6083, Facult\'e des Sciences, 37200 Tours, France}

\begin{abstract}
We study the behavior near $x_0$ of any positive
solution of (E) $-\Delta u=u^q$ in $\Omega$ which vanishes on
$\partial\Omega\setminus\{x_0\}$, where $\Omega\subset \RR^N$ is a
smooth domain, $q\geq (N+1)/(N-1)$ and
$x_0\in\partial\Omega$. Our results are based upon \textit{a priori}\/
estimates of solutions of (E) and existence, non-existence and
uniqueness results for solutions of some nonlinear elliptic equations
on the upper-half unit sphere. 
{\it To cite this article: M.-F. Bidaut-V\'eron, A.C. Ponce, L. V\'eron, C. R.
Acad. Sci. Paris, Ser. I XXX (2006).}

\vskip 0.5\baselineskip

\selectlanguage{francais}
\noindent{\bf R\'esum\'e}
\vskip 0.5\baselineskip
\noindent
Nous \'etudions le comportement
quand  $x$ tend vers $x_0$ de toute solution positive de (E) $-\Delta
 u=u^q$ dans $\Omega$ qui s'annule sur
 $\partial\Omega\setminus\{x_0\}$, o\`u $\Omega\subset \RR^N$ est un
 domaine r\'egulier, $q\geq  (N+1)/(N-1)$ et $x_0\in\partial\Omega$. Nos
 r\'esultats sont fond\'es sur des estimations a priori des
 solutions de (E), et des r\'esultats d'existence, de non existence
 et d'unicit\'e de solutions de certaines \'equations elliptiques
 non lin\'eaires sur la demi-sph\`ere unit\'e. 
 {\it Pour citer cet article~:  M.-F. Bidaut-V\'eron, A.C. Ponce,
   L. V\'eron, C. R. Acad. Sci. Paris, Ser. I XXX (2006).}

\end{abstract}
\end{frontmatter}

\selectlanguage{francais}
\section*{Version fran\c{c}aise abr\'eg\'ee}
Soit $\Omega$ un ouvert r\'egulier de $\RR^N$, $N\geq 4$, tel que $0
\in \partial\Omega$. \'Etant donn\'e $q > 1$, nous consid\'erons une fonction $u\in C^2(\Omega) \cap 
C(\overline\Omega\setminus\{0\})$ qui v\'erifie
\eqref{Eq1}. Nous nous int\'eressons \`a la description du comportement
de $u$ au voisinage de $0$.

Nous distinguerons les trois valeurs critiques de $q$ donn\'ees par
\eqref{q}. Si $1<q<q_1$, le comportement en $0$ des solutions
est d\'ecrit dans \cite{BV-Vi}; aussi supposerons-nous le plus souvent
$q\geq q_1$. 
Si $u$ est une solution de \eqref{Eq1} dans $\RR^N_+$ de la forme
$u(x)=u(r,\sigma)=r^{-2/(q-1)}\omega(\sigma)$, alors $\omega$
v\'erifie l'\'equation \eqref{En3}.
Dans ce cas, nous avons le r\'esultat suivant:

\smallskip
\begin{theoreme} 
$(i)$ Si $1< q\leq q_1$, le probl\`eme \eqref{Eq1} n'admet aucune solution.

$(ii)$ Si $q_1<q<q_3$, \eqref{Eq1} admet une unique solution,
  not\'ee $\omega_0$.

$(iii)$ Si $q\geq q_3$, \eqref{Eq1} n'admet aucune solution.
\end{theoreme}
\smallskip

Le r\'esultat d'unicit\'e d\'ecrit en $(ii)$ est en fait un cas
particulier d'un r\'esultat plus g\'en\'eral:

\smallskip
\begin{theoreme}
Pour tous $1 < q \leq q_3$ et $\lambda\in \RR$, il existe au plus une solution
positive de \eqref{En4}.
\end{theoreme}
\smallskip

Ce r\'esultat demeure si, dans \eqref{En4}, $S^{N-1}_+$ est remplac\'e
par une boule dans $\RR^N$, et $\Delta'$ par le laplacien ordinaire.

\smallskip

Par simplicit\'e, nous pouvons supposer que $\partial\RR_+^{N}$ est
l'hyperplan tangent \`a $\Omega$ en $0$. Le th\'eor\`eme ci-dessous
donne une classification des singularit\'es isol\'ees du probl\`eme \eqref{Eq1}:

\smallskip
\begin{theoreme}
Soit $q\geq q_1$, avec $q \neq q_2$. Supposons que la solution
$u$ du probl\`eme \eqref{Eq1} v\'erifie
\begin{equation}\label {F5}
0\leq u(x)\leq C\abs x^{-2/(q-1)}\quad\forall x\in \Omega\cap B_a(0),
\end {equation}  
pour $C, a> 0$. Si $q_1 \leq q < q_3$, ou bien $u$ est continue
en $0$, ou bien
\begin{equation}\label {F7}
u(r,\sigma) = 
\begin{cases}
r^{-(N-1)} \, \big(\log{(1/r)}\big)^{\frac{1-N}{2}} \big(k_N \sigma_1 + o(1) \big) &
\mbox{si } q = q_1, \\
r^{-2/(q-1)} \big( \omega_0(\sigma) + o(1) \big) & \mbox{si } q_1 <
 q < q_3,
\end{cases}
\end{equation}
lorsque $r \to 0$, uniform\'ement par rapport \`a $\sigma \in S^{N-1}_+$;
$k_N$ est une constante qui d\'epend seulement de $N$.\\
Si $q\geq q_3$, $u$ est continue en $0$.
\end{theoreme}
\smallskip

L'estimation a priori \eqref{F5} est obtenue pour $q_1 \leq q <q_2$:

\smallskip
\begin{theoreme}
Si $q_1\leq q< q_2$, toute solution $u$ de \eqref{Eq1} v\'erifie
\eqref{F5} pour $C=C(N,q,\Omega)>0$. 
\end{theoreme}
\smallskip

Les d\'emonstrations d\'etaill\'ees sont pr\'esent\'ees dans \cite{BVPV}.


\selectlanguage{english}

\section{Introduction and main results}
\label{sec1}

Let $\Omega$ be a smooth open subset of $\RR^N$, $N \geq 4$, such that
$0 \in \partial\Omega$ and let $q>1$. Assume that $u\in C^2(\Omega) \cap
C(\overline\Omega\setminus\{0\})$ is a solution of  
\begin{equation}\label{Eq1}
\left\{
\begin{alignedat}{2}
-\Delta u & = u^q && \quad\mbox{in } \Omega,\\
u & \geq 0 && \quad\mbox{in } \Omega,\\
u & = 0 && \quad\mbox{on } \partial\Omega\setminus\{0\}.
\end{alignedat}
\right.
\end{equation}  
Our goal in this paper is to describe the behavior of $u$ in a
neighborhood of $0$.

This problem has similar features with the case where $x_0 \in
\Omega$, which has been studied by Gidas-Spruck~\cite{GS}. In our case,
we encounter three critical values of $q$ in
describing the local behavior of $u$: 
\begin{equation}\label{q}
q_1:=\frac{N+1}{N-1}, \quad q_2:=\frac{N+2}{N-2} 
\quad \mbox{and} \quad q_3:=\frac{N+1}{N-3}.
\end {equation} 
When $1<q<q_1$, it is proved in \cite{BV-Vi} that for every solution
$u$ of \eqref{Eq1} there exists $\alpha \geq 0$ (depending on $N$ and $u$) such that 
\begin{equation}\label{En2}
u(x)=\alpha \abs x^{-N}\rho(x) \, \big(1+o(1) \big)\quad \mbox{as } x\to 0,
\end {equation}
where $\rho(x)=\dist (x,\partial\Omega)$, $\forall x \in \Omega$.
For this reason, we shall mainly restrict ourselves to $q\geq
q_1$. 

Let us first consider the case where $\Omega = \RR^N_+$ and we look
for solutions of \eqref{Eq1} of the form
$u(x)=u(r,\sigma)=r^{-2/(q-1)}\omega(\sigma)$, where $r = |x|$ and
$\sigma \in S_+^{N-1}$. An easy computation shows that $\omega$ must satisfy
\begin{equation}\label{En3}
\left\{
\begin{alignedat}{2}
-\Delta'\omega & =\ell_{N,q}\omega + \omega^q && \quad \mbox{in }S^{N-1}_+,\\
\omega & \geq 0 && \quad\mbox{in } S^{N-1}_+,\\
\omega & = 0 && \quad\mbox{on }\partial S^{N-1}_+,
\end{alignedat}
\right.
\end {equation}  
where $\Delta'$ denotes the Laplacian in $S^{N-1}$ and 
$\ell_{N,q}= \frac{2(N-q(N-2))}{(q-1)^2}$.
Concerning equation \eqref{En3}, we prove

\smallskip

\begin{theorem}\label{Th1} 
$(i)$ If $1<q\leq q_1$, then \eqref{En3} admits no
  positive solution.

$(ii)$ If $q_1<q< q_3$, then \eqref{En3} admits a unique
positive solution.

$(iii)$ If $q\geq q_3$, then \eqref{En3} admits no positive solution.
\end{theorem}
\smallskip

One of the main ingredients in the proof of Theorem~\ref{Th1} $(ii)$
is the following

\smallskip

\begin{theorem}\label{Th'1} 
If $1 < q \leq q_3$ and $\lambda\in \RR$, then there exists at most one positive solution of 
\begin{equation}\label{En4}
\left\{
\begin{alignedat}{2}
-\Delta'v & =\lambda v + v^q && \quad\mbox{in }S^{N-1}_+,\\
v & =0 && \quad\mbox{on }\partial S^{N-1}_+.
\end{alignedat}
\right.
\end{equation}  
\end{theorem}


We now return to the case where $\Omega \subset \RR^N$ is an arbitrary
smooth set such that $0 \in \partial\Omega$. For simplicity, we may
assume that $\partial\RR_+^N$ is the tangent 
hyperplane of $\Omega$ at $0$. Using Theorem~\ref{Th'1}, we provide
a classification of isolated singularities of solutions of \eqref{Eq1}:
\smallskip

\begin{theorem}\label{Th2} 
Let $q\geq q_1$, $q \neq q_2$, and let $u$ be a solution of
\eqref{Eq1}. Assume that $u$ satisfies 
\begin{equation}\label {En5}
0\leq u(x)\leq C\abs x^{-2/(q-1)}\quad\forall x\in \Omega\cap B_a(0),
\end {equation}  
for some $C, a> 0$. If $q_1 \leq q < q_3$, then either $u$ is
continuous at $0$ or
\begin{equation}\label {En7}
u(r,\sigma) = 
\begin{cases}
r^{-(N-1)} \, \big(\log{(1/r)}\big)^{\frac{1-N}{2}} \big(k_N \sigma_1 + o(1) \big) &
\mbox{if } q = q_1, \\
r^{-2/(q-1)} \big( \omega_0(\sigma) + o(1) \big) & \mbox{if } q_1 <
 q < q_3,
\end{cases}
\end{equation}
as $r \to 0$, uniformly with respect to $\sigma \in S^{N-1}_+$;
$k_N$ denotes a constant depending only on $N$ and $\omega_0$ is the
unique positive solution of \eqref{En3}.\\
If $q\geq q_3$, then $u$ is continuous at $0$.
\end{theorem}
\smallskip

\begin{remark} 
{\rm 
We do not know whether Theorem~\ref{Th2} is true when $q=q_2$. In this
case, the equation is conformally invariant and thus other techniques
are required. If $\Omega= \RR_+^N$, then it can be proved
that any solution of \eqref{Eq1} depends only on the 
  variables $r = |x|$ and $\theta=\cos^{-1}(x_1/\abs x)$.}
\end{remark}
\smallskip

The next result establishes the existence of an \textit{a priori}
estimate for the solutions of \eqref{Eq1}. According to
Theorem~\ref{Th3} below, assumption \eqref{En5} is always fulfilled when
$q_1 \leq q < q_2$:

\smallskip

\begin{theorem}\label{Th3} 
Let $ q_1\leq q<q_2$ and let $u$ be a solution of \eqref{Eq1}. Then,
\begin{equation}\label {En9}
0\leq u(x)\leq C\rho(x)\abs x^{-2/(q-1)-1}\quad \forall x \in \Omega \cap B_1(0),
\end {equation}
where $C$ depends on $N$, $q$ and $\Omega$.
\end{theorem}
\smallskip

\begin{remark}
{\rm
According to the Doob Theorem~\cite{D}, any positive superharmonic
function $v$ in 
$\Omega$ satisfies $\int_\Omega |\Delta v| \, \rho < \infty$ and admits a
boundary trace, which is a Radon measure on $\partial\Omega$. If $u$
is a solution of \eqref{Eq1}, then its trace must be of the form $k
\delta_{x_0}$, for some $k \geq 0$. We may have $k > 0$ if $1 < q <
q_1$ (see \cite{BV-R}), but $k$ is necessarily equal to $0$ if $q \geq
q_1$. Indeed, by the maximum principle, $u$ satisfies $u \geq k
P_\Omega(x,0)$, where $P_\Omega$ denotes the Poisson potential of
$\Omega$. Since $u^q \in 
L_\rho^1(\Omega)$ (by the Doob Theorem), we must have $k = 0$ if $q
\geq q_1$.}
\end{remark}
\smallskip

Detailed proofs will appear in \cite{BVPV}.


\section{Sketch of the proofs}

\noindent
\textit{Proof of Theorem~\ref{Th1}.}\/
Assertion $(i)$ is proved by multiplying \eqref{En3} by $\phi(\sigma) =
\sigma_1$. Note that $\phi$ is the first eigenfunction of $-\Delta'$ on
$S_+^{N-1}$, with eigenvalue $\lambda_1 = N-1$. Integrating the
resulting expression over $S^{N-1}_+$, and using the fact that $1 < q
\leq q_1\Longrightarrow\ell_{N,q}\geq\lambda_1$, we obtain $(i)$.\\
The existence part in $(ii)$ is obtained by using the Mountain Pass
Theorem; the uniqueness is a consequence of Theorem~\ref{Th'1}.\\
Assertion $(iii)$ can be deduced from the following Poho\v zaev-type
identity:

\smallskip
\begin{proposition}\label{Po} 
Assume $q>1$. Then, any solution of \eqref{En4} satisfies
\begin{multline*}
\frac{N-3}{q+1}\left(q-q_3\right)
\int_{S^{N-1}_+}\abs{\nabla' v}^2\phi\,d\sigma
-\frac{(N-1)(q-1)}{q+1}\left(\lambda+ \frac{N-1}{q-1}\right)\int_{S^{N-1}_+}{}
v^{2}\phi\,d\sigma= \\
= - \int_{\partial S^{N-1}_+}{}\abs{\nabla'v}^2 \,d\tau.
\end{multline*}
\end{proposition}

This identity is obtained by computing the divergence of the vector
 field $P=\langle\nabla'\phi,\nabla'v\rangle\nabla'v$, where $\nabla'$
 is the gradient on $S^{N-1}$, and then using the fact that the first
 eigenfunction satisfies $D^2\phi+\phi g_0=0$, where $g_0$ is the
 tensor of the standard metric on $S^{N-1}$. In order to establish $(iii)$, it
 suffices to observe that $\ell_{N,q} \leq - \frac{N-1}{q-1}
 \Longleftrightarrow q \geq q_3$.

\medskip 
\noindent
\textit{Proof of Theorem~\ref{Th'1}.}\/ We first notice that any
positive solution of (\ref{En4}) depends only on the variable
$\theta=\cos^{-1}(x_1/\abs x)\in [0,\pi/2]$; this follows from a
straightforward adaptation of the Gidas-Ni-Nirenberg moving plane method
to $S^{N-1}_+$ (see \cite{P}). Thus, $v$ satisfies 
 \begin{equation}\label{u1}
\left\{\begin{array} {l}
 v''+(N-2) \cot{\theta} \, v'+ \lambda v + v^q =0\quad\mbox{in }(0,\pi/2),\\
 v'(0)=0,\quad v(\pi/2)=0.
\end{array}\right. 
\end{equation}
Let  $w(\theta):= \sin^\alpha\theta \, v(\theta)$, where $\alpha>0$. By
 choosing $\alpha = 2(N-2)/(q+3)$, then $w$ satisfies 
\begin{equation}\label{u5}
\big(w'(\pi/2) \big)^2= \int_{0}^{\pi/2} G'(\theta) w^2(\theta) \, d\theta,
\end{equation}
where $G$ is a function of the form $G(\theta) = \sin^{\beta'}\!\theta \,
\big( \alpha_1 \sin^2{\theta} + \alpha_2 \big)$; the parameters
$\alpha_1,  \beta' \in \RR$  and $\alpha_2 \geq 0$ can be explicitly computed in
terms of $\lambda$, $N$ and $q$.\\
Assume, by contradiction, that $v_1$ and $v_2$ are two distinct
solutions of \eqref{u1}. Then,
\begin{equation}\label{u2}
\int_0^{\pi/2} v_1 v_2 \big(v_2^{q-1}-v_1^{q-1} \big) \, d\theta=0.
\end{equation}
Therefore, their graphs must intersect at some $\theta_0 \in (0, \pi/2)$. 
We claim that $v_1$ and $v_2$ intersect at least twice in $(0,
\pi/2)$. 
If there is only one intersection point, then assuming $v_2(0) > v_1(0)$ it can be shown
that there exists $\gamma \geq \left( \frac{w_2'(\pi/2)}{w_1'(\pi/2)} \right)^2$ such that the function $\theta \mapsto G'(\theta) \big(w_2^2(\theta) - \gamma w_1^2(\theta) \big)$ is nonnegative in $(0, \pi/2)$.
Thus, by \eqref{u5},
$$
0 < \int_0^{\pi/2} G'(\theta) \left[w^2_2(\theta) - \gamma w^2_1(\theta) \right] \, d\theta = \big( w_2'(\pi/2) \big) ^2 - \gamma \big( w_1'(\pi/2) \big) ^2 \leq 0.
$$
This is a contradiction. Therefore, $v_1$ and $v_2$ must intersect at least twice. This fact
leads to another contradiction by using the Shooting Method
(see \cite{KL}). Thus, $v_1 = v_2$ in $(0,\pi/2)$.

\smallskip 

\begin{remark}
{\rm
The method above follows the lines of the proof of
Kwong-Li~\cite{KL}. 
}
\end{remark}
\medskip

\noindent 
\textit{Proof of Theorem~\ref{Th2}.}\/
It follows from  methods developed in \cite {GS} and \cite{BV-V}. For
simplicity, we shall assume that $a = 1$ and 
 $\partial\Omega \cap B_1 = \partial\RR_+^N \cap B_1$. We set
 $$
w(t,\sigma) = r^{2/(q-1)} u(r,\sigma),\quad t=\log{(1/r)} \in 
 (0,\infty)\times S^{N-1}_+:= Q.
 $$
Then, $w$ satisfies
\begin{equation}\label{as1}
w_{tt}-\left(N-2\frac{q+1}{q-1}\right)w_t+\Delta'w+\ell_{N,q}w+w^q=0
\quad \mbox{in } Q
 \end {equation}
 and $w$ vanishes on $(0,\infty)\times \partial S^{N-1}_+$. Since $w$
 is uniformly bounded on $Q$, standard \textit{a priori}\/ estimates for
 elliptic problems yield
 $$
\abs{\partial^k_t\nabla'^jw}\leq M_{k,j}\quad \mbox{in } (1,\infty)\times S^{N-1}_+
 $$
for any integers $k,j \geq 0$, where $\nabla'^j$ stands for the covariant
derivative on $S^{N-1}$. Thus, the trajectory $\cT_w = \big\{w(t,\cdot):t\geq
1\}$ is relatively compact in 
$C^2\big(\overline{S^{N-1}_+} \big)$. 
Multiplying \eqref{as1} by $w_t$ and integrating over $S_+^{N-1}$, we obtain 
\begin{equation}\label{as2}
\frac{d}{dt}H(t)=\left(N-2\frac{q+1}{q-1}\right)
\int_{S^{N-1}_+}w_t^2 \, d\sigma,
\end {equation}
where
$$
H(t):=\frac{1}{2}\int_{S^{N-1}_+}
\left(w_t^2-\abs{\nabla' v}^2-\ell_{N,q}w^2+\frac{2}{q+1}w^{q+1}\right)d\sigma.
$$
Since $q \neq q_2$, we know that $N-2(q+1)/(q-1) \neq 0$. Thus,
iterated energy estimates imply that $w_t(t,\cdot),\,w_{tt}(t,\cdot)\to 0$
in $L^2(S^{N-1}_+)$ as $t\to\infty$. Therefore, the limit set
$\Gamma_w$ of $\cT$ is a connected subset of the set of
solutions of \eqref{En3}. By Theorem~\ref{Th1}, we deduce
that 
$$
\Gamma_w = 
\begin{cases}
\{0\}  & \mbox{if $q = q_1$ or $q \geq q_3$},\\
\{0\} \mbox{ or } \{\omega_0\} & \mbox{if $q_1 < q < q_3$.} 
\end{cases}
$$
Then, a linearization argument as in
\cite{BV-V} leads to the conclusion if $q > q_1$.\\
We now consider the case $q = q_1$; we borrow some ideas from
\cite{BV-R} and \cite{V}. We first prove, by ODE techniques, that  
\begin{equation}\label{as3}
X(t):=\int_{S^{N-1}_+} w(t,\cdot)\phi\, d\sigma \leq Ct^{-(N-1)/2}.
\end {equation}
Using \eqref{En5} and the boundary Harnack inequality (see
\cite{CFMS}), we derive  
\begin{equation}\label{as4}
 0\leq w(t,\sigma)\leq  Ct^{-(N-1)/2}\quad \mbox{in } (1,\infty)\times S^{N-1}_+.
\end {equation}
Set $\eta(t,\sigma):=t^{(N-1)/2}w(t,\sigma)$. We verify as above
that the limit set $\Gamma_\eta$ in $C^2\big(\overline{S^{N-1}_+}\big)$ of the
trajectory $\cT_\eta$ of $\eta$ is an interval of the form 
$\big\{\kappa\phi:0\leq \kappa_0\leq \kappa\leq
\kappa_1\big\}$. In order to show that $\cT_\eta$ is reduced to a
single point, we prove that $\norm{r(t,\cdot)}_{L^2}\leq Ct^{-1}$, where
$$
r(t, \cdot ):=\eta(t, \cdot )-z(t) \phi \quad \mbox{and} \quad z(t)=\int_{S^{N-1}_+}\eta(t,.)\phi\,d\sigma.
$$
Writing the equation satisfied by $z$ as a non-homogeneous second
order linear ODE, we prove that either $z(t)\to 0$, which implies
that $u$ is continuous at $0$, or $z(t)\to \tilde k_N$ as
$t\to\infty$, for some constant depending only on $N$.
\medskip 

\noindent
\textit{Proof of Theorem~\ref{Th3}.}\/ It is an application of the
Doubling Lemma Method introduced in \cite{PQS}, from which we derive
the following local estimate:

\smallskip
\begin{lemma}\label{lem} 
Let $1<q<q_2$ and let $u$ be a solution of \eqref{Eq1}. Then, for
every $x_0 \in \partial\Omega \setminus \{0\}$ and $0 < R < |x_0|$, we have
\begin{equation}
0 \leq u(x) \leq  C\big(R-\abs{x-x_0} \big)^{-2/(q-1)} \quad\forall x\in B_R(x_0) \cap \Omega,
\end{equation}
for some constant $C>0$ depending only on $\Omega$.
\end{lemma}
\smallskip

Apply this lemma with $x_0 \in \partial\Omega\setminus\{0\}$ and $R =
|x_0|/2$. Using elliptic regularity theory, we obtain
$$
0 \leq u(x)\leq C\rho(x)\abs x^{-2/(q-1) - 1} \quad \forall x\in\Omega \mbox{
  such that } 0 < \rho(x)<|x|/2.
$$
If $\rho(x)\geq |x|/2$, then we use Gidas-Spruck's internal estimates
(see \cite{GS}). We thus obtain \eqref{En9}.



\section*{Acknowledgements}
The authors are grateful to E.~N.~Dancer for his invaluable comments.


\begin{thebibliography}{99}


\bibitem{BV-R} M.-F. Bidaut-V\'eron, Th. Raoux, { Asymptotics of
  solutions of some nonlinear elliptic systems},
  Comm. Partial Differential Equations 21 (1986), 1035--1086.

\bibitem{BVPV} M.-F. Bidaut-V\'eron, A.C. Ponce, L. V\'eron, in
  preparation.

\bibitem{BV-V} M.-F. Bidaut-V\'eron, L. V\'eron, { Nonlinear
  elliptic equations on compact Riemannian manifolds and asymptotics
  of Emden equations}, Invent. Math. 106 (1991), 489--539. 

\bibitem{BV-Vi} M.-F. Bidaut-V\'eron, L. Vivier, { An elliptic
  semilinear equation with source term involving boundary measures: the
  subcritical case}, Rev. Mat. Iberoamericana 16 (2000), 477--513.

\bibitem{CFMS} L. Cafarelli, E. Fabes, S. Mortola, S. Salsa, {
 Boundary behavior of nonnegative 
 solutions of elliptic operators in divergence form},
 Indiana Univ. Math. J. 30 (1981), 621--640.
 
\bibitem{D} J. Doob, Classical potential theory and its probabilistic
  counterpart, Springer, London, 1984.

\bibitem{GS} B. Gidas, J. Spruck, { Global and local behavior of
  positive solutions of nonlinear elliptic equations}, Comm. Pure
  Appl. Math. 34 (1981), 525--598.


\bibitem{KL} M.K. Kwong, Y. Li, { Uniqueness of radial solutions of
  semilinear elliptic equations}, Trans. Amer. Math. Soc. 333 (1992),
  339--363.
 
\bibitem{P} P. Padilla, { Symmetry properties of positive solutions of
  elliptic equations on symmetric domains}, Appl. Anal. 64 (1997), 153--169.

\bibitem{PQS} P. Pol\'a\v cik, P. Quittner, Ph. Souplet, { Singularity
  and decay estimates in superlinear problems via Liouville type
  theorems. Part I: Elliptic equations and systems}, Duke Math. J., to
  appear. 

\bibitem{V} L. V\'eron, { Comportement asymptotique des solutions
  d'\'equations elliptiques semi-lin\'eaires dans $\RR^N$},
  Ann. Math. Pura Appl. 127 (1981), 25--50.






\end{thebibliography}
\end{document}